\documentclass[doublespacing]{bmcart}
\usepackage{amsmath}
\usepackage{amsfonts}
\usepackage{amssymb}
\usepackage[utf8]{inputenc} 

\def\includegraphics{}

\startlocaldefs

\DeclareMathOperator{\e}{e}

\newtheorem{lemma}{Lemma}
\newtheorem{theorem}{Theorem}

\newtheorem{definition}{Definition}
\newtheorem{remark}{Remark}
\newtheorem{Example}{Example}
\newenvironment{example}{\begin{Example}\normalfont}{\end{Example}}

\endlocaldefs

\begin{document}

\begin{frontmatter}

\begin{fmbox}
\dochead{Research}


\title{Numerical Algorithm for Nonlinear Delayed Differential Systems of $n$th Order}


\author[
   addressref={aff1},                   
   corref={aff1},                       
   email={josef.rebenda@ceitec.vutbr.cz}   
]{\inits{J}\fnm{Josef} \snm{Rebenda}}
\author[
   addressref={aff2},
   email={smarda@feec.vutbr.cz}
]{\inits{Z}\fnm{Zden\v{e}k} \snm{\v{S}marda}}


\address[id=aff1]{
  \orgname{Brno University of Technology, CEITEC BUT}, 
  \street{Purky\v{n}ova 656/123},                     %
  \postcode{61200}                                
  \city{Brno},                              
  \cny{Czech Republic},                                    
  {josef.rebenda@ceitec.vutbr.cz}
}
\address[id=aff2]{%
  \orgname{ Brno University of Technology, Faculty of Electrical Engineering and Communication},
  \street{Technick\'a 10},
  \postcode{61600}
  \city{Brno},
  \cny{Czech Republic},
  {smarda@feec.vutbr.cz}
}


\begin{artnotes}
\end{artnotes}

\end{fmbox}


\begin{abstractbox}

\begin{abstract} 
The purpose of this paper is to propose a semi-analytical technique convenient for numerical approximation of solutions of the initial value problem for $p$-dimensional delayed and neutral differential systems with constant, proportional and time varying delays.
The algorithm is based on combination of the method of steps and the differential transformation. Convergence analysis of the presented method is given as well.
Applicability of the presented approach is demonstrated in two examples: A system of pantograph type differential equations and a system of neutral functional differential equations with all three types of delays considered.  
Accuracy of the results is compared to results obtained by the Laplace decomposition algorithm, the residual power series method and Matlab package DDENSD. Comparison of computing time is done too, showing reliability and efficiency of the proposed technique.
\end{abstract}




\end{abstractbox}
%

\end{frontmatter}



\section{Introduction}
Systems of functional differential equations (FDEs), in particular delayed or neutral differential equations, are often used to model processes in the real world.
To give some examples, we mention models in population dynamics \cite{gyori1989}, neuromechanics \cite{insperger2015}, machine tool vibrations \cite{kalmar} etc. Further models and details can be found for instance in monographs \cite{hale} or \cite{kolmanovski}.

Semi-analytical methods expecting solutions to problems with delays in a series form have been studied in the last two decades. Methods such as the variational iteration method (VIM) \cite{chen}, the Adomian decomposition method (ADM) \cite{cocom}, the homotopy perturbation method (HPM) \cite{shakeri}, the homotopy analysis method (HAM) \cite{duarte2016} and also methods based on the Taylor theorem such as the differential transformation (DT) \cite{rebenda2017b}, the Taylor collocation method \cite{bellour} and the Taylor polynomial method \cite{sezer} have been developed to approximate solutions to different types of problems for FDEs. Other ways to use series approach in solving FDEs are e.g. the method of polynomial quasisolutions \cite{cherepennikov2010}, \cite{cherepennikov2013} or the functional analytic technique (FAT) \cite{petropoulou2007}, \cite{petropoulou2009}.

The main aim of the work is to apply the combination of the method of steps and DT as a convenient tool for finding an approximate solution to the initial value problem for functional differential systems used in dynamical models. Convergence analysis and error estimates of the method are investigated as well. We give some experimental results in section \ref{applications} to show that the algorithm produces reliable results with the same efficiency as or better efficiency than the reference methods.
 
\section{Methods}
The main idea of our approach is to combine the differential transformation and the general method of steps.


The differential transformation is an active research topic during the last years. As examples of recently published results, we mention research papers \cite{samajova}, \cite{rebenda2017a}, \cite{yang}, \cite{samajova2017} and \cite{rebenda2017c}. These papers as well as other publications contain new algorithms and their applications to solving different types of problems involving differential equations.

\begin{definition}
The differential transformation of a real function $u(t)$ at a point $t_0 \in \mathbb R$ is $\mathcal D \{ u(t) \} [t_0] = \{ U(k)[t_0] \}_{k=0}^{\infty}$,
where the $k$th component $U(k)[t_0]$ of the differential transformation of the function $u(t)$ at $t_0$ is defined as
\begin{equation}\label{c2}
U(k)[t_0] = \frac{1}{k!} \left[ \frac{d^ku(t)}{dt^k} \right]_{t=t_0},
\end{equation}
assuming that the original function $u(t)$ is analytic.
\end{definition}
\begin{definition}
 The inverse differential transformation of $\{ U(k)[t_0] \}_{k=0}^{\infty}$ at $t_0$ is defined as
\begin{equation}\label{c3}
u(t) = \mathcal D^{-1} \Bigl\{ \{ U(k)[t_0] \}_{k=0}^{\infty} \Bigr\} [t_0]= \sum_{k=0}^{\infty}U(k)[t_0] (t-t_0)^k. 
\end{equation}
\end{definition}
In applications, the function $u(t)$ is usually expressed in the form of finite series 
\begin{equation}\label{3f}
u(t) = \sum_{k=0}^{N}U(k)[t_0] (t-t_0)^k. 
\end{equation}
In section \ref{applications} we use the following transformation formulas which are derived from definitions \eqref{c2}, \eqref{c3} and proved in \cite{smardadiblikkhan2013}.
\begin{lemma}{}\label{lemma1}
Assume that $W(k)$,  $U(k)$ and $U_i (k)$ are the $k$th components of the differential transformations of functions $w(t)$, $u(t)$ and
$u_i (t)$, $i=1,2$, at $t_0 \in \mathbb{R}$, respectively, and let $q, q_j \in (0,1)$, $j=1,2$. Moreover, assume that $t_0=0$.
\begin{align*}
\text{If }& w(t) = {\displaystyle \frac{d^nu(t)}{dt^n}},   \text{ then }   W(k) = {\displaystyle \frac{(k+n)!}{k!}}U(k+n). \\[2mm]
\text{If }& w(t) = u_1 (t) u_2 (t),  \text{ then }  W(k) = \sum_{l=0}^k U_1 (l) U_2 (k-l).\\
\text{If }& w(t) =  u(qt),  \text{ then }  W(k) =q^k U(k).\\
\text{If }& w(t) =  u_1(q_1t)u_2(q_2t),  \text{ then } W(k) = \sum_{l=0}^k q_1^lq_2^{k-l}U_1(l)U_2(k-l).\\
\text{If }& w(t) \! = \! {\displaystyle \frac{d^mu(qt)}{d(qt)^m}} \! = \! \left. \frac{d^m u(t)}{dt^m} \right|_{t=qt} \! = u^{(m)}(qt), \text{then } \! W(k) \! = \! {\displaystyle \frac{(k+m)!}{k!}}q^k U(k+m).
\end{align*}
\begin{align*}
\text{If } w(t) &=  t^n,  \text{ then } W(k) =\delta (k-n), \text{ where }  \delta ( k - n ) = \delta_{kn} \text{ (Kronecker delta)}. \\[1mm]
\text{If } w(t) &=  \e^{\lambda t},  \text{ then }  W(k) = {\displaystyle \frac{\lambda^k}{k!}}.\\
\text{If } w(t) &= \cos t, \text{ then }
 C(k) = 
 \begin{cases}
(-1)^{\frac{k}{2}} \frac{1}{k!} & \text{if } k=2n, n \in \mathbb{N}_0 \\
  0 & \text{if } k=2n+1, n \in \mathbb{N}_0
 \end{cases}. \\
\text{If } w(t) &= \sin t,  \text{ then }
  S(k) = 
 \begin{cases}
(-1)^{\frac{k-1}{2}} \frac{1}{k!} & \text{if } k=2n+1, n \in \mathbb{N}_0 \\
  0 & \text{if } k=2n, n \in \mathbb{N}_0
 \end{cases}.
\end{align*}

\end{lemma} 

\begin{remark}{}
Transformation formulas for shifted arguments $w(t)=u(t-a)$ are often proved and applied in papers. However, using these formulas to solving initial value problems for delayed differential equations is not convenient since the uniqueness of solutions is violated. The reason is that the values of the initial vector function 
for $t < 0$ are not taken into account.
\end{remark}

One of the drawbacks of common approach to the differential transformation is that there is no use of direct transformation formulas for equations with nonlinear terms containing unknown function $u(t)$, e.g. $f(u)= \e^{\cos{u}}$ or $f(u) = \sqrt{1+u^4}$.

Fortunately, the corresponding transformations can be calculated using the Adomian polynomials $A_n$ in which each solution $u_i$ is replaced  by the corresponding components $U_i (k)$ of the differential transformation $\{ U_i (k) \}_{k=0}^{\infty}$, see \cite{smardakhan2015}. Suppose that $F(k)$ is the $k$th component of the differential transformation of a nonlinear term $f(u)$, then
\begin{align}\label{PP4}
F(k) &=  \sum_{n=0}^{\infty} A_n(U(0),U(1),\dots, U(n))\delta(k-n) = A_k(U(0),U(1),\dots, U(k)) \nonumber \\
&=  \frac{1}{k!} \frac{d^k}{dt^k}  \left[ f\left(  \sum_{l=0}^{\infty} U(l) t^l \right)\right]_{t=0}, \quad k \geq 0. 
\end{align}
Recently, it turned out that there is another possible way how to work with nonlinearities in DT \cite{rebenda2018c}.


The second method, the method of steps, enables us to replace the terms involving constant or time-dependent delays by initial vector function 
 and its derivatives. Then the original initial value problem for a system of delayed or neutral differential equations is simplified to the initial problem for a system of ordinary differential equations.
Details on the method of steps can be found for instance in monographs  \cite{hale}, \cite{kolmanovski} or \cite{bellen}.

\section{Results}

The subject of our interest is a system of $p$ functional differential equations of $n$th order with multiple delays $\alpha_1 (t), \dots, \alpha_r (t)$ in the following form:
\begin{align}\label{1}
{\mathbf u}^{(n)}(t) \! = &\mathbf{f}(t,\! \mathbf{u}(t),\! \mathbf{u}'(t),\dots,\! \mathbf{u}^{(n-1)}(t),\! \mathbf{u}_1(\alpha_1(t)), \! \mathbf{u}_2(\alpha_2(t)), \dots,\! \mathbf{u}_r(\alpha_r(t)) ),
\end{align}
where
$\mathbf{u}^{(n)} (t) =(u_1^{(n)}(t), \dots, u_p^{(n)}(t))^T\ , \mathbf{u}^{(k)} (t) = (u_1^{(k)}(t), \dots, u_p^{(k)}(t))$, $k=0,1, \ldots, n-1$ and $\mathbf{f} = (f_1, \dots, f_p)^T$ are $p$-dimensional vector functions,  $\mathbf{u}_i(\alpha_i(t))= (\mathbf{u}(\alpha_i(t)),\mathbf{u}'(\alpha_i(t)),\dots,\mathbf{u}^{(m_i)}(\alpha_i(t)))$ are  $(m_i \cdot p)$-dimensional vector functions, $m_i \leq n$, $i=1,2,\dots,r$, $r \in \mathbb N$ and $f_j \colon [0,\infty) \times \mathbb{R}^{np} \times \mathbb{R}^{\omega p}$ are continuous real functions for $j=1,2,\dots,p$, where $\omega = \sum\limits_{i=1}^r m_i$.

We consider three types of delays $\alpha_i$:

\begin{enumerate}
\item[1.]
$\alpha_i(t) = q_i t$, where $q_i \in (0,1)$ (proportional delay).
\item[2.]
 $\alpha_i(t) = t-\tau_i$, where $\tau_i>0$ is a real constant (constant delay).
\item[3.]
$\alpha_i(t)=t-\tau_i(t)$, where $\tau_i(t) \geq \tau_{i0}>0$ for $t>0$ is a real function (time-dependent or time-varying delay).
\end{enumerate}

Let $t^*= \min\limits_{1\leq i \leq r}\Bigl\{\inf\limits_{t>0} \bigl(\alpha_i(t)\bigr)\Bigr\} \leq 0$,
$m= \max\{m_1,m_2,\dots,m_r\} \leq n$.  In case $m=n$ system \eqref{1} is a neutral system, otherwise it is a delayed differential system.\\
If $t^*<0$, an initial vector function $\Phi(t) = (\phi_1(t), \dots,\phi_p(t))^T$ must be assigned to the system \eqref{1} on the interval $[t^*, 0]$.
Moreover, we assume that $\phi_j(t) \in C^n([t^*,0],\mathbb{R})$ for $j=1,\dots,p$.

We look for a solution of the system \eqref{1} with the following initial conditions
\begin{equation}\label{6}
\mathbf{u}(0)=\mathbf{v}_0,\mathbf{u}'(0)=\mathbf{v}_1, \dots, \mathbf{u}^{(n-1)}(0) = \mathbf{v}_{n-1}
\end{equation}
and the initial vector function $\Phi(t)$ on interval $[t^*, 0]$ satisfying
\begin{equation}\label{6'}
\Phi(0) = \mathbf{u}(0), \dots, \Phi^{(n-1)}(0)= \mathbf{u}^{(n-1)}(0).
\end{equation}
We solve initial value problem \eqref{1}, \eqref{6} and \eqref{6'} subject to the following hypotheses:
\begin{enumerate}
\item[(H1)]
The functions $f_j$, $j=1, \ldots, p$ are analytic in $[0,T^*] \times \mathbb{R}^{np} \times \mathbb{R}^{\omega p}$.
\item[(H2)]
The initial value problem \eqref{1}, \eqref{6} and \eqref{6'} has a unique solution on some interval $[0, T^*]$.
\end{enumerate}

\begin{remark}{}
\label{r1}
Hypothesis (H2) is valid for example if the delay functions $\alpha_i$ are Lipschitz continuous on $[0,T^*]$, the functions $\phi_j, \phi'_j, \ldots, \phi_j^{(n)}$ are Lipschitz continuous on $[t^*,0]$ and the functions $f_j$ are continuous with respect to $t$ on $[0,T^*]$ and Lipschitz continuous with respect to the rest of the variables on $\mathbb{R}^{np} \times \mathbb{R}^{\omega p}$. More details and other types of sufficient conditions for existence of a unique solution can be found in publications \cite{kolmanovski}, Subchapters 3.2 and 3.3,   or \cite{bellen}, Subchapter 2.2..
\end{remark}

We start with the method of steps. We substitute the initial vector function $\Phi(t)$ and its derivatives in all places where the unknown functions with constant or time-dependent delays and derivatives of that functions take place. This turns the delayed system \eqref{1} into a system of ordinary differential equations or differential equations with proportional delays in case that the system \eqref{1} contains proportional delays.

For example, if $\alpha_1 (t) = t-\tau_1$, $\alpha_2 (t) = t-\tau_2$, $\alpha_3 (t) = q_3 t$ and $\alpha_4 (t) = t-\tau_4(t)$, applying the method of steps changes \eqref{1} into the system
\begin{equation}\label{7}
\mathbf{u}^{(n)}(t) = \mathbf{f}(t,\mathbf{u}(t),\dots,\mathbf{u}^{(n-1)}(t),\mathbf{\Phi}_1 (t-\tau_1),\mathbf{\Phi}_2 (t-\tau_2), \mathbf{u}_3(q_3 t), \mathbf{\Phi}_4(t-\tau_4(t)) ),
\end{equation}
where
\begin{align*}
&\mathbf{\Phi}_i(t-\tau_i ) = (\Phi(t-\tau_i),\Phi'(t-\tau_i),\dots,\Phi^{(m_i)}(t-\tau_i)), \ i=1,2,\\
& \mathbf{u}_3 (q_3 t) = (\mathbf{u}(q_3 t), \mathbf{u}'(q_3 t),\dots, \mathbf{u}^{(m_3)} (q_3 t)), \\
&\mathbf{\Phi}_4(t-\tau_4(t) )= (\Phi(t-\tau_4(t)),\Phi'(t-\tau_4(t)), \dots,\Phi^{(m_4)}(t-\tau_4(t))),
\end{align*}
and $m_l \leq n$ for $l=1,2,3,4$. Then we transform the initial conditions \eqref{6}. Definition \eqref{c2} gives
$$
\mathbf{U}(k) = \frac{1}{k!} \mathbf{u}^{(k)} (0).
$$
After applying the differential transformation, initial value problem for a system of FDEs is reduced to a system of recurrence algebraic relations
\begin{align}\label{8}
\mathbf{U}(k+n) = \mathcal F \Bigl( k, \mathbf{U}(k),\mathbf{U}(k+1), \dots, \mathbf{U}(k+n-1) \Bigr).
\end{align}
Solving this recurrence and then using the inverse transformation \eqref{c3}, we get an approximate solution of
the system \eqref{1} in the series form
$$
\mathbf{u}(t) = \sum_{k=0}^{\infty} \mathbf{U}(k)t^k.
$$
If $t^*<0$, we denote $t_{\alpha_i} = \inf\{ t: \alpha_i (t)>0 \}$ and $t_{\alpha} =\min\limits_{1 \leq i \leq r} \{ t_{\alpha_i}: t_{\alpha_i} \neq 0 \}$. Then the approximate solution $\mathbf{u}(t)$ is valid on the intersection of its convergence interval and the interval $[0,T^*] \cap [0,  t_{\alpha}]$, whereas $\mathbf{u}(t) = \Phi(t)$ on the interval $[t^*,0]$.  If $t^*=0$, the approximate solution $\mathbf{u}(t)$ is valid on the intersection of its convergence interval with $[0,T^*]$.


Now we formulate and prove two theorems on convergence and an error estimate of the approximate solution to the studied problem obtained using differential transformation.

\begin{theorem}\label{T1}
Let the hypotheses (H1) and (H2) be valid and denote $\mathbf{F}_k(t) = \mathbf{U}(k)t^k$. If there exists  a constant $\delta$, $0<\delta < 1$, and $k_0 \in \mathbb{N}$
such that $|| \mathbf{F}_{k+1}(t)|| \leq \delta || \mathbf{F}_{k}(t)||$ for all $k \geq k_0$, then the series $ \sum_{k=0}^{\infty} \mathbf{F}_k(t)$ converges to a unique solution on the interval $J=[0,\gamma]$, $\gamma \leq T^*$.
\end{theorem}
\vspace*{2mm}
Proof. Denote $C^{n}(J)$ the Banach space of vector-valued functions\\
$\mathbf{h}(t) = (h_1(t),h_2(t),\dots, h_p(t))^T$ with continuous derivatives up to order $n$ and norm 
$$
||\mathbf{h}(t)|| = \max_{i=1, \dots, p} \max_{j=0, \dots n} \max_{t \in J}|h_i^{(j)} (t)|.
$$
Denote
$$
\mathbf{S}_l = \sum_{k=0}^l \mathbf{F}_k(t).
$$
Now it is sufficient to prove that the sequence $\left\{\mathbf{S}_l \right\} $ is a Cauchy  sequence  in the Banach  space $C^{n}(J)$.
Considering
$$
||\mathbf{S}_{l+1}-\mathbf{S}_l|| = ||\mathbf{F}_{l+1}(t)|| \leq \delta ||\mathbf{F}_{l}(t)|| \leq \dots \leq \delta^{l-n_0+1}||\mathbf{F}_{n_0}(t)||,
$$
then, for every $l,m \in \mathbb{N},\ l \geq m >n_0$, we get
\begin{align} \label{e}
||\mathbf{S}_l-\mathbf{S}_m|| &= || \sum_{j=m}^{l-1} ( \mathbf{S}_{j+1}-\mathbf{S}_j ) || \leq \sum_{j=m}^{l-1}||\mathbf{S}_{j+1}-\mathbf{S}_j|| \leq \sum_{j=m}^{l-1} \delta^{j-n_0+1}||\mathbf{F}_{n_0}(t)|| \nonumber \\
&= \delta^{m-n_0+1}(1+\delta +\delta^2 + \dots + \delta^{l-m-1} )||\mathbf{F}_{n_0}(t)|| \nonumber \\
& =\frac{1-\delta^{l-m}}{1-\delta} \delta^{m-n_0+1} ||\mathbf{F}_{n_0}(t)||.
\end{align}
Since  $0<\delta <1$ it follows 
$$
\lim_{l,m \rightarrow \infty} ||\mathbf{S}_l-\mathbf{S}_m|| = 0.
$$
Therefore, $\left\{\mathbf{S}_l \right\} $ is a Cauchy  sequence  in the Banach space $C^{n}(J)$ and the proof is complete.
\begin{theorem}\label{T2}
Suppose that the assumptions of Theorem \ref{T1} are valid.
Then for the truncated series $\sum_{k=0}^m\mathbf{F}_k(t)$ the following error estimate holds:
$$
|| \mathbf{u}(t) - \sum_{k=0}^m\mathbf{F}_k(t) || \leq \frac {1}{1-\delta} \delta ^{m-m_0+1}  \max_{i=1, \dots, p} \max_{j=0, \dots, n} \left| \frac{m_0!}{(m_0 -j)!} {U}_i(m_0) \gamma^{m_0-j} \right|.
$$
for any  $m_0 \geq 0$, $m \geq m_0$. 
\end{theorem}
Proof. Without loss of generality we can choose $m_0 \geq n$, where $n$ is the order of the system \eqref{1}. From inequality \eqref{e} we have
\begin{align}\label{e1}
||\mathbf{S}_l -\mathbf{S}_m &|| \leq \frac{1-\delta^{l-m}}{1-\delta} \delta^{m-m_0+1} ||\mathbf{F}_{m_0}(t)|| \nonumber \\
&= \frac{1-\delta^{l-m}}{1-\delta} \delta^{m-m_0+1} \max_{i=1, \dots, p} \max_{j=0, \dots, n} \left| \frac{m_0!}{(m_0 -j)!} {U}_i(m_0) \gamma^{m_0-j} \right|,
\end{align}
for  $ l \geq m \geq m_0$. From $0 <\delta <1$ it follows $ (1-\delta^{l-m})< 1$. Hence inequality \eqref{e1} can be reduced to
$$
||\mathbf{S}_l-\mathbf{S}_m|| \leq \frac {1}{1-\delta} \delta ^{m-m_0+1} \max_{i=1, \dots, p} \max_{j=0, \dots, n} \left| \frac{m_0!}{(m_0 -j)!} {U}_i(m_0) \gamma^{m_0-j} \right|.
$$
Here we use the fact that for $l \rightarrow \infty$, $\mathbf{S}_l \rightarrow \mathbf{u}(t)$ and the proof is complete.
\begin{remark}
Recent results on error estimates and convergence of Taylor series can be found e.g. in paper \cite{warne2006}.
\end{remark}

\section{Applications and Discussion}\label{applications}
As the first application, we have chosen the initial problem which has been solved in the paper \cite{widatalla2012} using the Laplace decomposition method (LDM) and in the paper \cite{komashynska2014} using the residual power series method (RPSM).

\begin{example}{}
We are looking for a solution of a $3$-dimensional system of pantograph equations
 \begin{align}
  u_1' (t) &= 2 u_2 \Bigl( \frac{t}{2} \Bigr) + u_3 (t) - t \cos \Bigl( \frac{t}{2} \Bigr), \notag \\
  u_2' (t) &= 1 - t \sin (t) -2 u_3^2 \Bigl( \frac{t}{2} \Bigr),  \label{p1}\\
  u_3' (t) &= u_2 (t) - u_1 (t) - t \cos (t). \notag
  \end{align}
subject to the initial conditions
\begin{equation}\label{p2}
u_1 (0) = -1,\quad u_2(0)=0, \quad u_3(0)= 0.
\end{equation}
Since the system \eqref{p1} contains proportional delays only, we do not have to use the method of steps.
Applying DT formulas in Lemma \ref{lemma1} to \eqref{p1} we get a system of recurrence relations
\begin{align}
  (k+1) U_1 (k+1) &= 2 \frac{1}{2^k} U_2 (k) + U_3 (k) - \frac{1}{2^{k-1}} C(k-1), \notag \\
  (k+1) U_2 (k+1) &=  \delta (k) - S(k-1) -2 \sum\limits_{l=0}^{k} \frac{1}{2^{k}} U_3 (l) U_3 (k-l),  \label{p3}\\
  (k+1) U_3 (k+1) &= U_2 (k) - U_1 (k) - C(k-1). \notag
\end{align}
From the initial conditions we have $U_1(0)=-1$, $U_2(0)=0$, $U_3(0)=0$. Solving the system \eqref{p3} we get
\begin{align*}
k=0:\ U_1(1) &= 2 U_2(0)+ U_3 (0)  =0, \\
U_2(1) &= \delta (0) - 2 (U_3 (0))^2=  1,\\
U_3(1) &= U_2(0) - U_1(0) =1, \\
k=1: \ U_1(2) &= \frac{1}{2}\left( 2 \frac{1}{2} U_2(1)+U_3(1) - C(0) \right)=\frac{1}{2}, \\
U_2(2) &= \frac{1}{2}\left( \delta (1) - 2 \frac{1}{2} (U_3 (0) U_3 (1) + U_3 (1) U_3 (0)) \right) =  0,\\
U_3(2) &= \frac{1}{2}\left( U_2(1) - U_1(1) - C(0) \right)= 0.
\end{align*}
For $k \geq 2$, we find
 $$ \begin{array}{lll}
U_1(3)=0, \quad & U_1(4)=-\frac{1}{4!}, \quad & U_1(5)=0, \ \dots \\[2mm]
U_2(3) =-\frac{1}{2}, \quad & U_2(4) = 0, \quad & U_2(5)= \frac{1}{4!}, \ \dots \\[2mm]
U_3(3)=  -\frac{1}{3!}, \quad & U_3(4)= 0, \quad & U_3(5)= \frac{1}{5!}, \ \dots 
\end{array} $$
Application of the inverse differential transformation \eqref{c3} gives a solution to \eqref{p1}, \eqref{p2} in the form
\begin{align*}
u_1(t) &= -1 + \frac{1}{2} t^2 - \frac{1}{4!} t^4 + \ldots = - \sum\limits_{k=0}^{N} (-1)^k \frac{t^{2k}}{(2k)!}, \\
u_2(t) &= t - \frac{1}{2} t^3 + \frac{1}{4!} t^5 - \ldots = \sum\limits_{k=0}^{N} (-1)^k \frac{t^{2k+1}}{(2k)!}, \\
u_3(t) &=   t - \frac{1}{3!} t^3 + \frac{1}{5!} t^5 - \ldots = \sum\limits_{k=0}^{N} (-1)^k \frac{t^{2k+1}}{(2k+1)!}.
\end{align*}
If $N \rightarrow \infty$, the series converge to the Taylor expansions of the closed form solutions 
 $$
 u_1(t)= -\cos t, \quad u_2(t) = t \cos t, \quad u_3(t) = \sin t.
 $$
Comparison of absolute errors of the presented  DT technique with LDM and RPSM for $N=2$ is done in Table~\ref{table1}, Table~\ref{table2} and Table~\ref{table3}. We see that DT and RPSM produce the same results which are close to the values of the closed form solutions, whereas LDM shows significant deviations. Similar results we obtain in comparison of computing times, see Table~\ref{table4},
Table~\ref{table5} and Table~\ref{table6}.

\begin{table}
\caption{Error analysis of $u_1$  on $[0,1]$.}
\centering
\label{table1}       
\begin{tabular}{|c|c|c|l|l|l|} 
\hline
           &exact solution  & DT       & abs. errors   & abs. errors & abs. errors  \\\hline\hline
  t     &$-\cos t$            &$u_1$&          DT&          LDM          & RPSM \\
\hline 
0.2      &-0.9800665&    -0.9800666&          1.0E - 7&     8.904E - 5&    1.0E - 7 \\
0.4      &-0.9210609&    -0.9210666&          5.7E - 6&     1.511E - 3&    5.7E -6\\
0.6      &-0.8253335&    -0.8254000&          6.65E - 5&    8.051E - 3&    6.65E - 5  \\
0.8      &-0.6967067&    -0.6970666&          3.599E - 4&   2.665E - 2&    3.599E - 4 \\
1.0      &-0.5403023&    -0.5416666&          1.3642E - 3&  6.766E - 2&    1.3642E - 3  \\\hline
\end{tabular}
\end{table}

\begin{table}
\caption{Error analysis of $u_2$  on $[0,1]$.}
\centering
\label{table2}       
\begin{tabular}{|c|c|c|l|l|l|} 
\hline
           &exact solution  & DT       & abs. errors   & abs. errors & abs. errors  \\\hline\hline
  t     &$t\cos t$            &$u_2$&          DT&          LDM          & RPSM \\
\hline  
0.2      &0.1960133&     0.1960133&          0.0&           5.496E - 6&    0.0\\
0.4      &0.3684243&     0.3684266&          2.3E - 6&      1.808E - 4&    2.3E - 6\\
0.6      &0.4952013&     0.4952400&          3.87E - 5&     1.408E - 3&    3.87E - 5  \\
0.8      &0.5573653&     0.5576533&          2.89E - 4&     6.069E - 3&    2.89E - 4 \\
1.0      &0.5403023&     0.5416666&          1.3643 - 3&    1.890E - 2&    1.3643E - 3  \\\hline
\end{tabular}
\end{table}

\begin{table}
\caption{Error analysis of $u_3$  on $[0,1]$.}
\centering
\label{table3}       
\begin{tabular}{|c|c|c|l|l|l|} 
\hline
           &exact solution  & DT       & abs. errors   &  abs. errors &  abs. errors  \\\hline\hline
  t     &$\sin t$            &$u_3$&          DT&          LDM          & RPSM \\
\hline  
0.2      &0.1986693&     0.1986693&          0.0&           6.4558E - 5&    0.0\\
0.4      &0.3894183&     0.3894186&          3.0E - 7&      9.9595E - 4&    3.0E - 7\\
0.6      &0.5646424&     0.5646480&          5.60E - 6&     4.8397E - 3&    5.60E - 6  \\
0.8      &0.7173561&     0.7173973&          4.12E - 5&     1.4613E - 2&    4.12E - 5 \\
1.0      &0.8414709&     0.8416666&          1.957 - 3&     3.3917E - 2&    1.957 - 3  \\\hline
\end{tabular}
\end{table}

\begin{table}
\caption{Comparison of computing time for  $u_1$.}
\centering
\label{table4}       
\begin{tabular}{|c|l|l|l|} 
 \hline
  t              &DT          &LDM          & RPSM \\\hline\hline
0.2                &6.3E-05&     8.7E-04&    6.3E-05 \\
0.4                &6.5E-05&     6.7E-04&    6.5E-05\\
0.6                &6.5E-05&     6.9E-04&    6.5E-05  \\
0.8                &6.4E-05&     7.2E-04E&   6.4E-05  \\
1.0                &6.6E-05&     8.8E-04E&   6.6E-05  \\\hline
\end{tabular}
\end{table}

\begin{table}
\caption{Comparison of computing time for  $u_2$.}
\centering
\label{table5}       
\begin{tabular}{|c|l|l|l|} 
 \hline
  t              &DT          &LDM          & RPSM \\\hline\hline
0.2                &6.6E-05&     6.7E-04&    6.6E-05 \\
0.4                &6.4E-05&     6.7E-04&    6.4E-05\\
0.6                &6.7E-05&     6.7E-04&    6.7E-05  \\
0.8                &6.5E-05&     8.5E-04E&   6.5E-05  \\
1.0                &6.5E-05&     8.3E-04E&   6.5E-05  \\\hline
\end{tabular}
\end{table}

\begin{table}
\caption{Comparison of computing time for  $u_3$.}
\centering
\label{table6}       
\begin{tabular}{|c|l|l|l|} 
 \hline
  t              &DT          &LDM          & RPSM \\\hline\hline
0.2                &6.6E-05&     7.7E-04&    6.6E-05 \\
0.4                &6.7E-05&     6.6E-04&    6.7E-05\\
0.6                &6.6E-05&     6.7E-04&    6.6E-05  \\
0.8                &6.6E-05&     8.3E-04E&   6.6E-05  \\
1.0                &6.7E-05&     8.5E-04E&   6.7E-05  \\\hline
\end{tabular}
\end{table}
 \end{example}

\begin{remark}{}
In paper \cite{widatalla2012}, using LDM authors obtained only approximate solutions of the initial value problem \eqref{p1}, \eqref{p2}. Applying RPSM, authors were able to find closed form solutions  in paper \cite{komashynska2014}. However, the calculations are too complicated and the residual functions (RPSM) and initial guesses (LDM) contain analytical forms of functions $\sin$ and $\cos$, which means that these methods are not convenient for use in a purely numerical software.
\end{remark}

As the second application, we have chosen a system with all three types of delays considered to show reliability and efficiency of the proposed approach in solving difficult tasks. 

\begin{example}{}\label{ex2}
Let us solve a nonlinear system of neutral delayed differential equations
\begin{align}\label{p8}
u'''_1 &= u_1'''(t-2)u_1\left(\frac{t}{3}\right)+\sqrt[3]{u_1^2}+ u'_2\Big(t-\frac{1}{2} \e^{-t}\Big), \nonumber \\
u'''_2 &= \frac{1}{2}u_2'''\left(\frac{t}{2}\right)  +u'_2(t-1)u_1\left(\frac{t}{3}\right)
\end{align}
 with initial functions
\begin{align}\label{P3}
 \phi_1(t) &= \e^t,\nonumber \\
 \phi_2(t) &=t^2
 \end{align}
 for $ t \in [-2,0]$, and initial conditions
\begin{align}\label{P4}
u_1(0)&=  1, \ u'_1(0)= 1,\ u_1''(0)=1,\nonumber \\
u_2(0)&= 0,\ u'_2(0)= 0,\ u''_2(0)=2.
\end{align}
Following the method of steps we get
\begin{align}\label{p9}
u'''_1 &= \e^{(t-2)}u_1\left(\frac{t}{3}\right)+ +\sqrt[3]{u_1^2}+ 2t- \e^{-t}, \nonumber\\
u'''_2 &= \frac{1}{2}u_2'''\left(\frac{t}{2}\right) +2(t-1)u_1\left(\frac{t}{3}\right).
\end{align}
The system \eqref{p9} cannot be solved by classical DT approach because of the nonlinear term $f(u) = \sqrt[3]{u_1^2}$, hence
we apply modified Adomian formula for differential transformation components to the nonlinear term $f(u)$. 
Applying DT to \eqref{p9} we get the system
\begin{align}
 (k+1)(k+2)(k+3)U_1(k+3) &= \e^{-2}\sum_{l=0}^k \frac{1}{l!}\left(\frac{1}{3}\right)^{k-l}U_1(k-l) +F_1(k) \notag \\
 &+2\delta (k) -\frac{(-1)^k}{k!}, \label{p10}\\
(k+1)(k+2)(k+3)\left(1 -\frac{1}{2^{k+1}}\right) &U_2(k+3) = 2 \sum_{l=0}^k \delta(l-1) \left(\frac{1}{3}\right)^{k-l}U_1(k-l) \notag \\
&- \frac{2}{3^k}U_1(k), \label{p11}
\end{align}
where $F_1(k)$ is the $k$th component of the transformed function $f(u) = \sqrt[3]{u_1^2}$. Applying the formula \eqref{PP4} and the transformed initial conditions 
\begin{align*}
U_1(0)& =1,\  U_1(1)=1,\ U_1(2)=\frac{1}{2}, \\
U_2(0)& =0,\  U_2(1)=0,\ U_2(2)=1,
\end{align*}
 we obtain
\begin{align*}
F_1(0) &= \sqrt[3]{U_1^2(0)} = 1, \\
F_1(1) &= \frac{2}{3} \frac {U_1(1)}{\sqrt[3]{U_1(0)}} = \frac{2}{3},\\
F_1(2) &= \frac{2}{3} \frac {U_1(2)}{\sqrt[3]{U_1(0)}} -\frac{1}{9} \frac {U_1^2(1)}{\sqrt[3]{U_1^4(0)}}= \frac{2}{9},\\
\vdots
\end{align*}
 Solving the system of recurrence relations \eqref{p10}, \eqref{p11} we get
\begin{align*}
k=0:\ U_1(3) &= \frac{1}{6}\Big( \e^{-2}U_1(0) +F_1(0) + 1\Big)= \frac{2+\e^{-2}}{6},\\
U_2(3) &= \frac{1}{3}\Big( -2U_1(0)\Big)=-\frac{2}{3},\\
k=1:\ U_1(4) &= \frac{\e^{-2}}{24}\Big( \frac{1}{3}U_1(1) + U_1(0)\Big)+\frac{1}{24} F_1(1) + 1\Big) = \frac{4\e^{-2}+5}{72},\\ 
U_2(4) &= \frac{1}{18}\Big(2U_1(0)-\frac{2}{3}U_1(1)\Big)=\frac{2}{27},\\
k=2:\ U_1(5) &= \! \frac{\e^{-2}}{60}\Big( \frac{1}{9}U_1(2) \! + \! \frac{1}{3}U_1(1) \!+ \! \frac{1}{2}U_1(0)\Big) \!+ \! \frac{1}{60}\Big( \! F_1(2) \! - \! \frac{1}{2} \Big) \! = \! \frac{16\e^{-2}-5}{1080},\\
U_2(5) &= \frac{2}{105}\Big(\frac{2}{3}U_1(1) -\frac{2}{9}U_1(2) \Big) =\frac{2}{189}.
\end{align*}
Applying the inverse differential transformation we obtain an approximate solution to the initial value problem \eqref{p8}, \eqref{P3}, \eqref{P4}:

\begin{align*}
u_1(t) &= 1+t +\frac{1}{2} t^2+ \frac{2+\e^{-2}}{6}t^3+ \frac{4\e^{-2}+5}{72}t^4+  \frac{16\e^{-2}-5}{1080}t^5 +\dots,\\
u_2(t) &= t^2- \frac{2}{3}t^3 + \frac{2}{27}t^4 + \frac{2}{189}t^5 + \dots.
\end{align*}

As we do not know exact solution of the given problem, we are limited to compare approximate solutions. Comparison of values obtained by the proposed approach and values obtained by Matlab package DDENSD in Table~\ref{table7} shows a good correspondence between the results.
Comparing computing times in Table~\ref{table8}, we can see that the presented method produces reliable results much faster than Matlab package DDENSD.

\begin{table}
\caption{Comparison of values of solution components obtained by DT and Matlab.}
\centering
\label{table7}       
\begin{tabular}{|c|c|c|c|c|}
\hline
Method    &   DT         & DT       & Matlab    & Matlab  \\\hline\hline
 t         &$u_1$         &$u_2$&          $u_1$&     $u_2$  \\
\hline
0.00      &1.0000&       0.0000&          1.0000&     0.0000 \\
0.05      &1.0513&       0.0024&          1.0513&     0.0024   \\
0.10      &1.1051&       0.0093&          1.1050&     0.0093   \\
0.15      &1.1618&       0.0203&          1.1614&     0.0203   \\
0.20      &1.2209&       0.0348&          1.2204&     0.0348   \\
0.25      &1.2832&       0.0524&          1.2822&     0.0524   \\
0.30      &1.3481&       0.0726&          1.3469&     0.0726  \\
0.35      &1.4160&       0.0951&          1.4146&     0.0951  \\\hline
\end{tabular}
\end{table}

\begin{table}
\caption{Comparison of computing time.}
\centering
\label{table8}       
\begin{tabular}{|c|c|c|c|c|}
\hline
Method    &   DT         & DT       & Matlab    & Matlab  \\\hline\hline
 t         &$u_1$         &$u_2$&          $u_1$&     $u_2$  \\
\hline
0.05      &7.9E-05&       6.5E-05&          6.2E-02&     6.2E-02   \\
0.10      &6.9E-05&       7.0E-05&          8.1E-02&     6.2E-02   \\
0.15      &7.1E-05&       6.8E-05&          6.2E-02&     6.2E-02   \\
0.20      &7.0E-05&       6.9E-05&          6.3E-02&     6.2E-02   \\
0.25      &7.0E-05&       7.0E-05&          6.2E-02&     6.2E-02   \\
0.30      &6.3E-05&       6.9E-05&          6.3E-02&     6.3E-02  \\
0.35      &6.9E-05&       6.8E-05&          6.2E-02&     6.3E-02  \\\hline
\end{tabular}
\end{table}
\end{example}

\begin{remark}{}
The system \eqref{p8} contains all three types of delay which were considered in this paper. Moreover, it contains a term which is nonlinear (nonpolynomial) in the dependent variable $u_1$. In this sense, the present paper contains more complicated systems in applications than  papers about other semi-analytical methods like VIM \cite{chen}, ADM \cite{cocom} or HPM \cite{shakeri}.
\end{remark}

\section{Conclusions}

The approach presented in this paper is an effective semi-analytical technique convenient for numerical approximation of a unique solution to the initial value problem for systems of functional differential equations, in particular delayed and neutral differential equations. Comparison of results was done against the Laplace decomposition method, the residual power series method and Matlab package DDENSD. 
The need of computational work is reduced compared to the other methods.
The differential transformation algorithm gives an approximate solution which is in a good concordance with reference results produced by Matlab. Under certain circumstances, it is possible to identify the unique solution to the initial value problem in the closed form. Further steps can be done in development of the presented technique for systems with distributed and state dependent delays. 


\begin{backmatter}

\section*{Competing interests}
  The authors declare that they have no competing interests.

\section*{Funding}
  The first author was supported by the Grant CEITEC 2020 (LQ1601) with financial support from the Ministry of Education, Youth and Sports of the Czech Republic under the National Sustainability Programme II. The work of the second author was supported by the Grant FEKT-S-17-4225 of Faculty of Electrical Engineering and Communication, Brno University of Technology.

\section*{Author's contributions}
 All authors contributed equally to the writing of this paper.


\bibliographystyle{bmc-mathphys} 
\bibliography{ADE_RebendaSmarda_2018.bib} 





\end{backmatter}
\end{document}